# CONJECTURES ON PARTITIONS OF INTEGERS AS SUMMATIONS OF PRIMES

Florentin Smarandache, Ph D
Associate Professor
Chair of Department of Math & Sciences
University of New Mexico
200 College Road
Gallup, NM 87301, USA
E-mail: smarand@unm.edu

Abstract.
In this short note many conjectures on partitions of integers as summations of prime numbers are presented, which are extension of Goldbach conjecture.

A) Any odd integer $n$ can be expressed as a combination of three primes as follows:

1) As a sum of two primes minus another prime: $n = p + q - r$, where $p, q, r$ are all prime numbers.
   Do not include the trivial solution: $p = p + q - q$ when $p, q$ are prime.
   For example:
   $$1 = 3 + 5 - 7 = 5 + 7 - 11 = 7 + 11 - 17 = 11 + 13 - 23 = ...;$$
   $$3 = 5 + 5 - 7 = 7 + 19 - 23 = 17 + 23 - 37 = ...;$$
   $$5 = 3 + 13 - 11 = ...;$$
   $$7 = 11 + 13 - 17 = ...$$
   $$9 = 5 + 7 - 3 = ...;$$
   $$11 = 7 + 17 - 13 = ...;$$

   a) Is this a conjecture equivalent to Goldbach's Conjecture (any odd integer $\geq 9$ is the sum of three primes)?
   b) Is the conjecture true when all three prime numbers are different?
   c) In how many ways can each odd integer be expresses as above?

2) As a prime minus another prime and minus again another prime: $n = p - q - r$, where $p, q, r$ are all prime numbers.
   For example:
   $$1 = 13 - 5 - 7 = 17 - 5 - 11 = 19 - 5 - 13 = ...;$$
   $$3 = 13 - 3 - 7 = 23 - 7 - 13 = ...;$$
   $$5 = 13 - 3 - 5 = ...;$$
   $$7 = 17 - 3 - 7 = ...;$$
   $$9 = 17 - 3 - 5 = ...;$$

$$11 = 19 - 3 - 5 = \ldots \ .$$

a) In this conjecture equivalent to Goldbach's Conjecture?
b) Is the conjecture true when all three prime numbers are different?
c) In how many ways can each odd integer be expressed as above?

B) Any odd integer $n$ can be expressed as a combination of five primes as follows:
3) $n = p + q + r + t - u$, where $p, q, r, t, u$ are all prime numbers, and $t \neq u$.
For example:
$$1 = 3 + 3 + 3 + 5 - 13 = 3 + 5 + 5 + 17 - 29 = \ldots;$$
$$3 = 3 + 5 + 11 + 13 - 29 = \ldots;$$
$$5 = 3 + 7 + 11 + 13 - 29 = \ldots;$$
$$7 = 5 + 7 + 11 + 13 - 29 = \ldots;$$
$$9 = 5 + 7 + 11 + 13 - 29 = \ldots$$
$$11 = 5 + 7 + 11 + 17 - 29 = \ldots \ .$$

a) Is the conjecture true when all five prime numbers are different?
b) In how many ways can each odd integer be expressed as above?

4) $n = p + q + r - t - u$, where $p, q, r, t, u$ are all prime numbers, and $t, u \neq p, q, r$.
For example:
$$1 = 3 + 7 + 17 - 13 - 13 = 3 + 7 + 23 - 13 - 19 = \ldots \ ;$$
$$3 = 5 + 7 + 17 - 13 - 13 = \ldots;$$
$$5 = 7 + 7 + 17 - 13 - 13 = \ldots;$$
$$7 = 5 + 11 + 17 - 13 - 13 = \ldots;$$
$$9 = 7 + 11 + 17 - 13 - 13 = \ldots;$$
$$11 = 7 + 11 + 19 - 13 - 13 = \ldots \ .$$

a) Is the conjecture true when all five prime numbers are different?
b) In how many ways can each odd integer be expressed as above?

5) $n = p + q - r - t - u$, where $p, q, r, t, u$ are all prime numbers, and $r, t, u \neq p, q$
For example:
$$1 = 11 + 13 - 3 - 3 - 17 = \ldots;$$
$$3 = 13 + 13 - 3 - 3 - 17 = \ldots;$$
$$5 = 5 + 29 - 5 - 5 - 17 = \ldots;$$
$$7 = 3 + 31 - 5 - 5 - 17 = \ldots;$$
$$9 = 3 + 37 - 7 - 7 - 17 = \ldots \ ;$$
$$11 = 5 + 37 - 7 - 7 - 17 = \ldots \ .$$

a) Is the conjecture true when all five prime numbers are different?
b) In how many ways can each odd integer be expressed as above?

6) $n = p - q - r - t - u$, where $p, q, r, t, u$ are all prime numbers, and $q, r, t, u \neq p$.

For example:
$$1 = 13 - 3 - 3 - 3 - 3 = ...;$$
$$3 = 17 - 3 - 3 - 3 - 5 = ...;$$
$$5 = 19 - 3 - 3 - 3 - 5 = ...;$$
$$7 = 23 - 3 - 3 - 5 - 5 = ...;$$
$$9 = 29 - 3 - 5 - 5 - 7 = ...;$$
$$11 = 31 - 3 - 5 - 5 - 7 = ... .$$

  a) Is the conjecture true when all five prime numbers are different?
  b) In how many ways can each odd integer be expressed as above?

GENERAL CONJECTURE:

Let $k \geq 3$, and $1 < s < k$ be integers. Then:

i) If $k$ is odd, any odd integer can be expressed as a sum of $k - s$ primes (first set) minus a sum of $s$ primes (second set) [such that the primes of the first set is different from the primes of the second set].

  a) Is the conjecture true when all $k$ prime numbers are different?
  b) In how many ways can each odd integer be expressed as above?

ii) If $k$ is even, any even integer can be expressed as a sum of $k - s$ primes (first set) minus a sum of $s$ primes (second set) [such that the primes of the first set is different from the primes of the second set].

  a) Is the conjecture true when all $k$ prime numbers are different?
  b) In how many ways can each even integer be expressed as above?

**REFERENCE**

[1]  Smarandache, Florentin, "Collected Papers", Vol. II, Moldova State University Press at Kishinev, article "Prime Conjecture", p. 190, 1997.